\documentclass[11pt]{amsart}

\usepackage{amsmath,amssymb}

\def\COMMENT#1{}
\def\TASK#1{}

\def\noproof{{\unskip\nobreak\hfill\penalty50\hskip2em\hbox{}\nobreak\hfill%
        $\square$\parfillskip=0pt\finalhyphendemerits=0\par}\goodbreak}
\def\endproof{\noproof\bigskip}
\newdimen\margin   
\def\textno#1&#2\par{%
    \margin=\hsize
    \advance\margin by -4\parindent
           \setbox1=\hbox{\sl#1}%
    \ifdim\wd1 < \margin
       $$\box1\eqno#2$$%
    \else
       \bigbreak
       \hbox to \hsize{\indent$\vcenter{\advance\hsize by -3\parindent
       \sl\noindent#1}\hfil#2$}%
       \bigbreak
    \fi}

\newtheorem{firstthm}{Proposition}
\newtheorem{thm}[firstthm]{Theorem}

\newtheorem{lemma}[firstthm]{Lemma}
\newtheorem{cor}[firstthm]{Corollary}

\newtheorem{defin}[firstthm]{Definition}
\newtheorem{conj}[firstthm]{Conjecture}

\begin{document}
\title{Approximate Hamilton decompositions of random graphs}
\date{\today}
\author{Fiachra Knox, Daniela K\"uhn and Deryk Osthus}
\thanks {D.~K\"uhn was supported by the EPSRC, grant no.~EP/F008406/1
and by the ERC, grant no.~258345. D.~Osthus was supported by the EPSRC, grant no.~EP/F008406/1.} 

\begin{abstract}
We show that if $pn \gg \log n$ the binomial random graph $G_{n,p}$ has an approximate Hamilton decomposition.
More precisely, we show that in this range $G_{n,p}$ contains a set of edge-disjoint Hamilton cycles 
covering almost all of its edges.
This is best possible in the sense that the condition that $pn \gg \log n$ is necessary.
\end{abstract}
\maketitle

\section{Introduction}

Consider the random graph model where one starts with an empty graph and
successively adds edges which are chosen uniformly at random.
One of the most striking results in the theory of random graphs is the fact that this random graph
acquires a Hamilton cycle as soon as its minimum degree is at least 2.
This was proved by Bollob\'as~\cite{Bol84} %
\COMMENT{also claimed but apparently not proven by Komlos and Szemeredi}
and then soon afterwards generalized by Bollob\'as and Frieze~\cite{BF85} who showed that
for any fixed $k$ the above random graph contains $k$ edge-disjoint Hamilton cycles as soon as its minimum degree is at least $2k$.
(If the minimum degree is odd, one can guarantee an additional perfect matching, this is also the case for several of the
conjectures and results discussed below.)

More recently, Frieze and Krivelevich~\cite{FK05} conjectured that a similar result even holds for arbitrary edge densities:
\begin{conj}[Frieze and Krivelevich \cite{FK05}] \label{mainconj}
For any $p$, whp the binomial random graph $G_{n,p}$ contains $\lfloor \delta(G_{n,p}) \rfloor /2$ edge-disjoint Hamilton cycles.
\end{conj}
Here we say that a property of a random graph on $n$ vertices holds \emph{whp} if the probability that it holds tends to 1 as $n$ tends to infinity.
The result from~\cite{BF85} implies that Conjecture~\ref{mainconj} holds if 
$pn \le \log n+ O(\log \log n)$. Frieze and Krivelevich~\cite{FK08} extended this to all $p$ with 
$pn=(1+o(1)) \log n$. They proved also an approximate version of Conjecture~\ref{mainconj} for constant edge probabilities.
\begin{thm}[Frieze and Krivelevich \cite{FK05}]\label{dense}
Let $0 < p < 1$ be constant. Then whp $G_{n, p}$ contains $(1-o(1))np/2$ edge-disjoint Hamilton cycles.
\end{thm}
As for constant $p$, whp all degrees are close to $np$, the result implies that $G_{n,p}$ can `almost'
be decomposed into edge-disjoint Hamilton cycles. 
As remarked in~\cite{FK10}, the proof of~\cite{FK05} also works as long as $p$ is a little larger than $n^{-1/8}$.
In this paper, we extend this result to essentially the entire range of $p$.
\begin{thm}\label{AppHamDecomp}
For any $\eta > 0$, there exists a constant $C$ such that if $p \geq \frac{C \log n}{n}$, 
then whp $G_{n, p}$ contains $(1 - \eta)np/2$ edge-disjoint Hamilton cycles.
\end{thm}
While finalizing this paper, we learned that this result was proven independently by Krivelevich (personal communication).
Theorem~\ref{AppHamDecomp} is best possible in the  sense that if we relax the condition on $p$, 
then $G_{n,p}$ can no longer be `almost' decomposed into Hamilton cycles. Indeed,  
suppose that $pn=C \log n$ for some fixed $C$. Then 
there is an $\varepsilon>0$ so that
whp the minimum degree of $G_{n,p}$ is at most
$(1-\varepsilon)np$ (see Exercise~3.4 in~\cite{Bollobasbook}). %

A version of Conjecture~\ref{mainconj} for random regular graphs of bounded degree
was proved by Kim and Wormald~\cite{KW01}:
if $r \ge 4$ is fixed, then whp a random $r$-regular graph contains $\lfloor r/2 \rfloor$ edge-disjoint
Hamilton cycles.\COMMENT{case $r$ odd is discussed within the paper}

Hypergraph versions of Theorem~\ref{dense} for Hamilton $\ell$-cycles
were also recently considered by Frieze and Krivelevich ~\cite{FK10} as well as
Frieze, Krivelevich and Loh~\cite{FKL10}.
A Hamilton $\ell$-cycle in this case consists of a cyclic ordering of the vertices together with a cyclic sequence of edges, 
each consisting of $k$ consecutive vertices,  such that consectutive edges intersect in exactly 
$\ell$ vertices.

A  `deterministic' version of these results was recently proved by 
Christofides, K\"uhn and Osthus~\cite{CKO}: Suppose that $G$ is a $d$-regular graph on $n$ vertices, 
where $d \ge (1+\varepsilon) n/2$ and $n$ is sufficiently large. Then $G$ contains $(1-\varepsilon)d/2$ edge-disjoint Hamilton cycles. This approximately confirms a conjecture of 
Nash-Williams~\cite{Nash-Williams71b} which states that any $d$-regular graph where $d \geq n/2$ has $\lfloor d/2 \rfloor$ 
edge-disjoint Hamilton cycles.

A related conjecture of Erd\H{o}s (see~\cite{tom1})
states that almost all tournaments $G$
contain at least $\min \{ \delta^+(G),\delta^-(G) \}$ edge-disjoint Hamilton cycles.
It follows from results in~\cite{FK05} that this is approximately true.
K\"uhn, Osthus and Treglown~\cite{KOT} recently proved that we do not even require
$G$ to be random for this to hold: any almost regular tournament $G$ contains a set
of edge-disjoint Hamilton cycles which cover almost all edges of $G$.

The proof of Theorem~\ref{dense} in~\cite{FK05} actually works for any quasi-random graph. 
(Here a graph $G$  is defined to be quasi-random if it is almost regular and if the density of any large induced subgraph of $G$
is close to that of $G$.)
In contrast, our argument only seems to work for $G_{n,p}$ (for example, the proof of Lemma ~\ref{num2factors} breaks down for a quasi-random graph). It would be interesting to 
know whether our result can also be extended to some natural class of (sparse) quasi-random graphs.

\section{Notation and organization of the paper}

We consider the binomial random graph $G_{n,p_0}$ and let $$w_0 = \frac{np_0}{\log n}.$$
So using the notation and assumptions of Theorem \ref{AppHamDecomp}, $w_0 \geq C$. Our results will hold provided that $w_0$ is sufficiently large depending on $\eta$, which we will assume throughout.

The graph $G_1$ is generated by including each edge of $G_0$ independently at random with probability $(1-\frac{\eta}{4})$. We define $G_2 := G_0 \backslash G_1$. Note that $G_1 \sim G_{n, p_1}$ and $G_2 \sim G_{n, p_2}$, where $p_1 = (1 - \frac{\eta}{4})p_0$ and $p_2 = \frac{\eta p_0}{4}$. (Of course, the distributions of these random graphs are not independent of each other.)

The outline of the proof of Theorem ~\ref{AppHamDecomp} is as follows: We begin in Section 3
by stating and applying some large deviation bounds on the number of edges in certain subgraphs of $G_{n, p}$, which we will use later on. Then in Section 4, we will use Tutte's $r$-factor theorem to show that one may find a regular subgraph of $G_1$ whose degree is close to $np_1$ (the average degree of $G_1$). In Section 5 we show that this subgraph can almost be decomposed into 2-factors in such a way that each 2-factor has relatively few cycles. Finally in Section 6 we convert each of these 2-factors into Hamilton cycles, using the edges of $G_2$ (along with any edges of $G_1$ which were not included in our collection of 2-factors). This is achieved using an appropriate variant of the well-known rotation-extension technique, first introduced by P\'osa \cite{Posa}. The fact that each of the 2-factors originally had few cycles will allow us to place an upper bound on the number of edges needed to perform the conversions, and thus to show that the process can be completed before all of the edges of $G_2$ have been used up.

Throughout the paper we use the following notation: for a graph $G$ and sets $A, B$ of vertices of $G$, $e_G(A, B)$ is the number of edges of $G$ with one endpoint in $A$ and the other in $B$. Let $e_G(A) = e_G(A, A)$. On the other hand if $G$ is a graph then $e(G)$ denotes the number of edges, and for a graph $G$ with a spanning subgraph $H$, $G \backslash H$ denotes the graph obtained by removing the edges of $H$ from $G$. We omit floor and ceiling symbols in arguments where they do not have a significant effect. $\log$ denotes the natural logarithm.
\section{Large Deviation Bounds}

We will need the following Chernoff bounds, which are proved e.g., in Janson, {\L}uczak and Ruci\'nski \cite{sj-tl-ar}:
\begin{lemma} \label{chernoff}
Let $X \sim Bin(n, p)$. Then the following properties hold:
\begin{itemize}
\item[{\rm (i)}] If $\varepsilon < \frac{1}{2}$, then $\mathbb{P}(|X - np| \geq \varepsilon np) \leq e^{-\frac{\varepsilon^2 np}{3}}$.
\item[{\rm (ii)}] If $t \geq 7np$, then $\mathbb{P}(X \geq t) \leq e^{-t}$.
\end{itemize}
\end{lemma}

We can use these bounds to deduce some facts about the number of edges between subsets of vertices of a random graph, as follows:

\begin{lemma} \label{ABedges}
Let $G \sim G(n,p)$. Then whp, for any disjoint $A, B \subseteq [n]$, the following properties hold: Let $a = |A|$ and $b = |B|$. Then
\begin{itemize}
\item[{\rm (i)}] If $\left(\frac{1}{a} + \frac{1}{b}\right) \frac{\log n}{p} \geq \frac{7}{2}$, then $e_{G}(A, B) \leq 2(a+b)\log n$, and
\item[{\rm (ii)}] If $\left(\frac{1}{a} + \frac{1}{b}\right) \frac{\log n}{p} \leq \frac{7}{2}$, then $e_{G}(A, B) \leq 7abp$.
\end{itemize}
\end{lemma}

\proof ~(i) Let $X = e_{G}(A, B)$ and let $t = 2(a+b)\log n$. Since $X \sim Bin(ab, p)$, we have that $t \geq 7abp = 7\mathbb{E}X$. If $a+b < 3$ then the result is trivial; otherwise, by Lemma ~\ref{chernoff} we have that $\mathbb{P}(X \geq t) \leq e^{-t} = \left(\frac{1}{n^a n^b}\right)^2 \leq \frac{1}{n^3}\left(\frac{1}{n^a n^b}\right)$, and a union bound immediately gives the result.

~(ii) Similarly, we have $\mathbb{P}(X \geq 7abp) \leq e^{-7abp} \leq e^{-t}$ and the result follows.
\endproof

In an exactly similar way\COMMENT{(i) Let $X = e_{G}(A)$ and let $t = 2a\log n$. Since $X \sim Bin(\frac{a^2}{2}, p)$, we have that $t \geq \frac{7a^2p}{2} = 7\mathbb{E}X$. If $a < 2$ then the result is trivial; otherwise, by Lemma ~\ref{chernoff} we have that $\mathbb{P}(X \geq t) \leq e^{-t} = \left(\frac{1}{n^a}\right)^2 \leq \frac{1}{n^2}\left(\frac{1}{n^a}\right)$, and a union bound immediately gives the result.\\
(ii) Similarly, we have $\mathbb{P}(X \geq \frac{7a^2p}{2}) \leq e^{-\frac{7a^2p}{2}} \leq e^{-t}$ and the result follows.}, we can show that
\begin{lemma} \label{Aedges}
Let $G \sim G_{n, p}$. Then whp, for every $A \subseteq [n]$ the following properties hold: Let $a = |A|$. Then
\begin{itemize}
\item[{\rm (i)}] If $\frac{\log n}{ap} \geq \frac{7}{4}$, then $e_{G}(A) \leq 2a \log n$, and
\item[{\rm (ii)}] If $\frac{\log n}{ap} \leq \frac{7}{4}$, then $e_{G}(A) \leq \frac{7a^2p}{2}$.
\end{itemize}
\end{lemma}

For larger sets, Lemma \ref{ABlinear} gives a more precise result.
Note that
we allow $\alpha, \beta \rightarrow 0$ in the statement.
\begin{lemma} \label{ABlinear}
Let $G \sim G_{n, p}$. Then whp, for all pairs $A, B \subseteq [n]$ of disjoint sets the following property holds: Let $\alpha = |A|/n$ and $\beta = |B|/n$, and suppose that $\alpha \beta np \geq 700$. Then
$$\frac{13}{14}\alpha \beta n^2 p \leq e_{G}(A, B) \leq \frac{15}{14}\alpha \beta n^2 p.$$
\end{lemma}
\proof $e_G(A, B) \sim Bin(\alpha \beta n^2, p)$, so by Lemma ~\ref{chernoff},
$$\mathbb{P}\left(e_{G}(A, B) < \frac{13}{14}\alpha \beta n^2 p\right) \leq e^{-\frac{\alpha \beta n^2 p}{3 \cdot 14^2}} \leq e^{-\frac{700 n}{588}} \leq e^{-n \log (3.1)} = \frac{1}{(3.1)^n}.$$
A union bound over all $3^n$ possibilities now gives the result. 
The right-hand inequality follows in an exactly similar manner using the opposite tail estimate. \endproof

\section{Regular subgraphs of a random graph}
We first show that $G_1$ contains a regular subgraph of degree at least 
$$r_1 = np_0\left(1-\frac{3\eta}{4}\right),$$
where $r_1$ is taken to be even. 

We do this by using a theorem of Tutte: Let $G$ be an arbitrary graph, $r$ a positive integer and suppose that $S, T, U$ is a partition of $V(G)$. Then define 
$$R_r(S, T) = \sum_{v \in T} d(v) - e_{G}(S, T) + r(|S| - |T|)$$
and let $Q_r(S, T)$ be the number of \textit{odd} components of $G[U]$, where a component $C$ is odd if and only if $r|C| + e_{G}(C, T)$ is odd. (In our case it will often suffice to bound $Q_r(S, T)$ simply by the total number of components of $G[U]$.)
\begin{thm}[Tutte \cite{Tutte52}] \label{Tutte}
Let $r$ be a positive integer. A graph $G$ contains an $r$-factor if and only if $R_r(S, T) \geq Q_r(S, T)$ for every partition $S, T, U$ of $V(G)$.
\end{thm}
In order to apply Theorem ~\ref{Tutte} we will need an upper bound on $Q_r(S, T)$. We do this by observing that if $G[U]$ has many components, then it must contain a large \textit{isolated set} of vertices; that is, a set $A \subseteq U$ such that $e_G(A, U \backslash A) = 0$. This becomes useful when looking at a random graph, since (as we will prove in Lemma ~\ref{boundary}) it follows that whp $A$ has many neighbours in $S \cup T$. This gives a lower bound on $|S \cup T|$ in terms of $Q_r(S, T)$, and thus gives an upper bound on $Q_r(S, T)$ in terms of $|S \cup T| = |S| + |T|$.
\begin{lemma} \label{isoset}
Let $G$ be a graph with $v$ components. Then for any $v' \leq \frac{v}{2}$, there exists a set $W \subseteq V(G)$ which is isolated in $G$, such that $v' \leq |W| \leq \max \{ 2v', \frac{2|G|}{v} \}$.
\end{lemma}
\proof Call a component $C$ of $G$ \textit{small} if its order is at most $v'$, and \textit{large} otherwise. Suppose first that the union of all small components of $G$ also has order at most $v'$. Then the number of small components is at most $v'$, and hence there are at least $v - v' \geq \frac{v}{2}$ large components. So one of these components must have order at most $\frac{2|G|}{v}$, and we can set $W$ to be this component.

On the other hand, if the sum of the orders of small components is greater than $v'$ then we can form $W$ by starting with $\emptyset$ and adding small components one by one until $|W| \geq v'$. Now since the last component added has size at most $v'$, we have that $|W| \leq 2v'$.
\endproof

Given a graph $G$, define the \textit{boundary} $B_{G}(A)$ of a set $A \subseteq V(G)$ to be the set of vertices which are adjacent (in $G$) to some vertex of $A$, but are not themselves elements of $A$. We will use the following two lemmas to give a lower bound on $|B_{G_1}(A)|$:

\begin{lemma} \label{mindegree}
Whp,
\begin{itemize}
\item[{\rm (i)}] $\delta(G_1) \geq (1-\frac{\eta}{2})np_0$,
\item[{\rm (ii)}] $\Delta(G_1) \leq np_0$, and
\item[{\rm (iii)}] $\delta(G_2) \geq \frac{\eta np_0}{5}$.
\end{itemize}
\end{lemma}

\proof Note that for a vertex $x$ of $G_1$, $d(x) \sim Bin(n-1, p_1)$ and $\mathbb{E}(d(x)) =  (1-\frac{\eta}{4})(n-1)p_0$. By Lemma ~\ref{chernoff}, we have 
$$\mathbb{P}\left(d(x) \leq (1-\frac{\eta}{2})np_0\right) \leq e^{-(\frac{\eta}{5})^2 \frac{np_0}{3}} = n^{-\frac{\eta^2 w_0}{75}} \leq \frac{1}{n^2},$$ 
and a union bound gives the result. The bound on the maximum degree follows similarly, as does that on the minimum degree of $G_2$.
\endproof

\begin{lemma} \label{boundary}
The following holds whp: Let $H$ be a spanning subgraph of $G_0$, and let $A \subseteq [n]$ be nonempty. Let $\delta_A = \min_{x \in A} d_H(x)$. Then setting $a=|A|$ and $b = |B_H(A)|$, the following properties hold:
\begin{itemize}
\item[{\rm (i)}] If $\frac{\log n}{ap_0} \geq \frac{7}{2}$, then $b \geq \frac{a(\delta_A - 6 \log n)}{2 \log n}$. In particular, if $H = G_1$, then $b \geq a$.
\item[{\rm (ii)}] If $\frac{\log n}{ap_0} \leq \frac{7}{2}$, then $3a + b \geq \frac{\delta_A}{7p_0}$. In particular, if $H = G_1$, then $3a + b \geq \frac{n}{14}$.
\end{itemize}
\end{lemma}

\proof Let $B = B_H(A)$.

~(i) Then 
\begin{equation}
a\delta_A \leq \sum_{x \in A} d_H(x) = e_H(A, B) + 2e_H(A) \leq e_{G_0}(A, B) + 2e_{G_0}(A)
\label{bdryapprox}
\end{equation}
which by Lemmas ~\ref{ABedges}(i) and ~\ref{Aedges}(i) is at most $2(a+b) \log n + 4a \log n$. The final part follows from Lemma ~\ref{mindegree}(i).

~(ii) We claim that $e_{G_0}(A, B) \leq 7ap_0(a+b)$. Indeed, if $\left(\frac{1}{a} + \frac{1}{b}\right) \frac{\log n}{p_0} \leq \frac{7}{2}$, then by Lemma ~\ref{ABedges}(ii), $e_{G_0}(A, B) \leq 7abp_0 \leq 7ap_0(a+b)$. On the other hand, if $\left(\frac{1}{a} + \frac{1}{b}\right) \frac{\log n}{p_0} \geq \frac{7}{2}$, then $e_{G_0}(A, B) \leq 2(a+b) \log n \leq 7ap_0(a+b)$. Similarly\COMMENT{If $\frac{\log n}{ap_0} \leq \frac{7}{4}$, then by Lemma ~\ref{Aedges}(ii), $e_{G_0}(A) \leq 7 a^2 p_0$. On the other hand, if $\frac{\log n}{ap_0} \geq \frac{7}{4}$, then $e_{G_0}(A) \leq 2a \log n \leq 7 a^2 p_0$.}, Lemma ~\ref{Aedges} implies that $e_{G_0}(A) \leq 7a^2p_0$. Now 
$$a \delta_A \stackrel{(\ref{bdryapprox})}{\leq} e_{G_0}(A, B) + 2e_{G_0}(A) \leq 7ap_0(a+b) + 14a^2 p_0$$ and the result follows immediately. Again the final part follows by Lemma ~\ref{mindegree}(i).
\endproof
We will later use the above lemmas to show that taking successive neighbourhoods of a set will give us a set of size linear in $n$ in a reasonably short time. For now, they allow us to give a bound on the size of $Q_r(S, T)$ in terms of $|S|$ and $|T|$.

\begin{lemma} \label{qst}
In the graph $G_1$, whp, for any partition $S, T, U$ of $[n]$, $Q_{r_1}(S, T) \leq 150(|S| + |T|)$.
\end{lemma}
\proof Let $v$ be the number of components of $G_1[U]$. Consider first the case when $150 \leq v \leq \frac{n}{75}$. Then by applying Lemma ~\ref{isoset} to the graph $G_1[U]$, we have that there exists a set $W \subseteq [n]$, isolated in $G_1[U]$, such that $\frac{v}{2} \leq |W| \leq \frac{n}{75}$. Now by Lemma ~\ref{boundary} with $A = W$, 
$$|B_{G_1}(W)| \geq \min\{|W|, \frac{n}{14} - 3|W|\} \geq |W| \geq \frac{v}{2}.$$ 
But if $W$ is isolated in $G_1[U]$, then its boundary in $G_1$ lies entirely in $S \cup T$. So $\frac{v}{2} \leq |B_{G_1}(W)| \leq |S| + |T|$, and hence $Q_{r_1}(S, T) \leq v \leq 2(|S| + |T|)$.

Now consider the case $v \geq \frac{n}{75}$. Setting $v' = \frac{n}{150}$ in Lemma ~\ref{isoset}, we have a set $W$, isolated in $G_1[U]$, such that $\frac{n}{150} \leq |W| \leq \frac{n}{75}$. Again the boundary of $W$ in $G_1$ has size at least $|W| \geq \frac{n}{150}$, and hence $\frac{n}{150} \leq |S| + |T|$. So $150(|S| + |T|) \geq n$ and the result holds trivially, since $Q_{r_1}(S, T)$ cannot be greater than $n$.

Finally if $v \leq 150$, then $Q_{r_1}(S, T) \leq 150 \leq 150(|S| + |T|)$ unless we are in the trivial case $|S| = |T| = 0$. But if $S, T$ are both empty then $U = [n]$ and $G_1[U] = G_1$, which whp has only one component (this follows from the fact that a random graph $G_{n,p}$ with
 $pn \ge 2\log n$ is whp connected, as proved in~\cite{ER59}). Since $r_1$ is even  this component cannot be odd, whence $Q_{r_1}(S, T) = 0$. 
\endproof

\begin{lemma}\label{Tutte1}
In the graph $G_1$, whp, we have that $R_{r_1}(S, T) \geq Q_{r_1}(S, T)$ for any partition $S, T, U$ of $[n]$.
\end{lemma}
\proof Let $d_S, d_T$ be the average degrees of the vertices in $S, T$ respectively. Let $\rho = \frac{|T|}{|S|}$, and $s = |S|$. We consider the following cases:

\textbf{Case 1:} $\rho \leq \frac{1}{2}$. Then since $e_{G_1}(S, T) \leq d_T |T|$ and $|S| \geq 2|T|$, we have 
$$R_{r_1}(S, T) \geq r_1 (|S| - |T|) \geq \frac{r_1}{3}(|S| + |T|)$$ 
and for sufficiently large $n$, $\frac{r_1}{3}(|S| + |T|) \geq 150(|S| + |T|) \geq Q_{r_1}(S, T)$.

\textbf{Case 2:} $\rho \geq 4$. Observe that by the definition of $r_1$ and by Lemma ~\ref{mindegree}, we have
\begin{equation} \label{case2}
d_T - r_1 \geq \left(1 - \frac{\eta}{2}\right)np_0 - \left(1 - \frac{3\eta}{4}\right)np_0 = \frac{\eta np_0}{4}
\end{equation}
and
$$d_S - r_1 \leq np_0 - \left(1 - \frac{3\eta}{4}\right)np_0 = \frac{3\eta np_0}{4}.$$
Now since $e_{G_1}(S, T) \leq d_S |S|$ and $|T| \geq 4|S|$,
\begin{align*}
R_{r_1}(S, T) &\geq d_T|T| - d_S|S| + r_1 (|S| - |T|) = (d_T - r_1)|T| - (d_S - r_1)|S| \\
&\geq \frac{\eta n p_0}{4}(|T| - 3|S|) \geq \frac{\eta n p_0}{20}(|S| + |T|)
\end{align*}
which again is at least $Q_{r_1}(S, T)$ for sufficiently large $n$.

\textbf{Case 3:} $\frac{1}{2} \leq \rho \leq 4$ and $\left(\frac{1}{s} + \frac{1}{\rho s}\right)\frac{\log n}{p_1} \geq \frac{7}{2}$. In this case by Lemma ~\ref{ABedges} we have that $e_{G_1}(S, T) \leq 2 (\rho+1) s \log n$, and so it suffices to prove that
\begin{equation}
\rho s(d_T - 2\log n -r_1 -150) + s(r_1 - 2 \log n -150) \geq 0.
\label{TCase3}
\end{equation}
(\ref{TCase3}) holds if $d_T - 2\log n -r_1 -150 \geq 0$ and $r_1 - 2 \log n -150 \geq 0$. But the latter inequality holds since $r_1 = \left(1 - \frac{3\eta}{4}\right)np_0 = \left(1 - \frac{3\eta}{4}\right) w_0 \log n$, and the former since 
$$
d_T - r_1 \stackrel{(\ref{case2})}{\ge}  \frac{\eta np_0}{4} = \frac{\eta w_0 \log n}{4}  \geq 3 \log n,
$$
as $w_0 \geq \frac{12}{\eta}$.

\textbf{Case 4:} $\frac{1}{2} \leq \rho \leq 4$ and $(\frac{1}{s} + \frac{1}{\rho s})\frac{\log n}{p_1} \leq \frac{7}{2}$ and $\rho s \leq \frac{n}{30}$. In this case by Lemma ~\ref{ABedges} we have that $e_{G_1}(S, T) \leq 7\rho s^2\left(1-\frac{\eta}{4}\right)p_0$, and so it suffices to prove that $$\rho s(d_T -r_1 -150) + s(r_1 - 7\rho s\left(1-\frac{\eta}{4}\right)p_0 -150) \geq 0,$$
and hence it suffices that $d_T -r_1 -150 \geq 0$ and $r_1 - 7\rho s(1-\frac{\eta}{4})p_0 -150 \geq 0$. But the former inequality holds as before, and the latter since 
$$r_1 - 150 \geq \frac{14}{15}r_1 = \frac{14}{15}\left(1 - \frac{3\eta}{4}\right)np_0 \geq 28\rho s\left(1 - \frac{3\eta}{4}\right)p_0 
\geq 7\rho s\left(1-\frac{\eta}{4}\right)p_0.$$

\textbf{Case 5:} $\frac{1}{2} \leq \rho \leq 4$ and $\rho s \geq \frac{n}{30}$. Note that we still have $\frac{s}{n} \leq \frac{1}{\rho+1}$, since $S, T$ are disjoint. Now by Lemma ~\ref{ABlinear}, $$e_{G_1}(S, T) \leq \frac{15}{14}\rho s^2 \left(1-\frac{\eta}{4}\right)p_0 \leq \frac{15}{14}\frac{\rho}{\rho+1} sn \left(1-\frac{\eta}{4}\right)p_0 \leq \frac{6}{7}sn \left(1-\frac{\eta}{4}\right)p_0.$$ So it suffices to prove that $$\rho sn\left(1-\frac{\eta}{2}\right)p_0 - \frac{6}{7}sn \left(1-\frac{\eta}{4}\right)p_0 + (1 - \rho)sn\left(1-\frac{3\eta}{4}\right)p_0 - 150(\rho+1)s \geq 0,$$ i.e., that $\frac{\eta \rho}{4} + (1-\frac{3\eta}{4}) - \frac{6}{7} (1-\frac{\eta}{4}) - \frac{150(\rho+1)}{np_0} \geq 0$, which is true if $\eta$ is not too large (which we can assume without loss of generality).
\endproof

\begin{cor}\label{regular}
Whp, $G_1$ contains an even-regular subgraph of degree $r_1 = (1 - \frac{3\eta}{4})np_0$.
\end{cor}
\proof This follows immediately from Theorem ~\ref{Tutte} and Lemma ~\ref{Tutte1}.
\endproof


\section{2-factors of regular subgraphs of a random graph}
In this section we will show that any even-regular subgraph of $G_0$ of sufficiently large degree contains a 2-factor with fewer than $\frac{\kappa n}{\log n}$ cycles, where 
\begin{equation}
\kappa = 2\log \left(\frac{16}{\eta}\right).  \label{kappa}
\end{equation}
It will follow immediately that we can decompose almost all of our regular subgraph into 2-factors with at most this many cycles. Roughly, our strategy will be to show that the number of 2-factors with many cycles in the original graph is rather small; smaller, in fact, than the minimum number of 2-factors which an even-regular graph of degree $r_1$ must contain.

\begin{lemma}\label{2factor} Let $k_0 = \frac{\kappa n}{\log n}$. Then whp, for any $r$-regular subgraph $H \subseteq G_0$ with $r \geq 2np_0 e^{-\frac{\kappa}{2}}$, $H$ contains a $2$-factor with at most $k_0$ cycles.
\end{lemma}
To prove Lemma ~\ref{2factor} we will need a number of further lemmas. We use Lemmas ~\ref{fracsum} and ~\ref{num2factors} to bound the number of 2-factors in $G_{n, p}$ with many cycles, while Lemma ~\ref{perm} gives a bound on the total number of 2-factors in $H$.
\begin{lemma}\label{fracsum}
For any $k$ and for $n \geq 3k$, we have 
$$\sum \prod_{i=1}^k \frac{1}{a_i} \leq \frac{k}{n}(\log n)^{k-1},$$
where the sum is taken over all ordered $k$-tuples $(a_1, a_2, \ldots, a_k)$ such that $a_1 + \ldots + a_k = n$ and $a_i \geq 3$ for each $i \in [n]$.
\end{lemma}
\proof We proceed by induction on $k$. The case $k = 1$ is trivial since both sides equal $\frac{1}{n}$. Supposing that the result holds for $k-1$, we have
$$\sum \prod_{i=1}^k \frac{1}{a_i} =  \sum_{a_k = 3}^{n - 3(k-1)}  \frac{1}{a_k} \sum \prod_{i=1}^{k-1} \frac{1}{a_i},$$
where again the second sum on the right-hand side is taken over all ordered $(k-1)$-tuples 
$(a_1, \ldots, a_{k-1})$ such that $a_1 + \ldots + a_{k-1} = n - a_k$ and $a_i \geq 3$ for all $i \in [k-1]$. By induction, this is bounded above by
\begin{align*}
&\sum_{a_k = 3}^{n-3(k-1)}  \frac{1}{a_k} \frac{k-1}{n-a_k}(\log (n-a_k))^{k-2} \\
&= \frac{k-1}{n} \sum_{a_k = 3}^{n-3(k-1)} \left(\frac{1}{a_k} + \frac{1}{n-a_k}\right)(\log (n-a_k))^{k-2} \\
&\leq \frac{k-1}{n}\left((\log n)^{k-2} \left(\sum_{a_k = 3}^{n-3} \frac{1}{a_k}\right) + \sum_{a_k = 3}^{n-3} \frac{1}{n-a_k} (\log (n-a_k))^{k-2}\right) \\
&\leq \frac{k-1}{n}\left((\log n)^{k-1} + \frac{1}{k-1}(\log n)^{k-1}\right) = \frac{k}{n} (\log n)^{k-1}, \\
\end{align*}
where the last inequality follows from the fact that $\log n = \int_1^n \frac{1}{x} \, dx$ and $\frac{1}{k-1}(\log n)^{k-1}$ $= \int_1^n \frac{1}{x} (\log x)^{k-2} \, dx$.
\endproof

\begin{lemma}\label{num2factors} Let $G \sim G_{n, p}$. Then whp, for any $k \geq \log n$ the number $A_k$ of $2$-factors in $G$ with at least $k$ cycles satisfies
$$A_{k+1} < \frac{n! (\log n)^{2k} p^n}{k! 2^k}.$$
\end{lemma}
\proof Note that it suffices to show that if $A'_k$ is the number of 2-factors in $G$ with \textit{exactly} $k$ cycles, then 
\begin{equation}
\mathbb{E}(A'_{k}) \leq \frac{(n-1)! (\log n)^{k-1} p^n}{(k-1)! 2^{k}}. \label{EAk1}
\end{equation}
Indeed, we then have 
\begin{align*}
\mathbb{E}(A_{k+1}) &= \sum_{i=k+1}^{\frac{n}{3}} \mathbb{E}(A'_{i}) 
\leq \sum_{i=k+1}^n \frac{(n-1)! (\log n)^{i-1} p^n}{(i-1)! 2^{i}} \\
&\leq \sum_{i=k+1}^n \frac{(n-1)! (\log n)^k p^n}{k! 2^{i}} 
\leq \frac{(n-1)! (\log n)^k p^n}{k! 2^k} \\
\end{align*}
and Markov's inequality implies that 
$$
\mathbb{P}\left(A_{k+1} \geq \frac{n! (\log n)^{2k} p^n}{k! 2^k}\right) \leq \frac{1}{n (\log n)^{k}} \leq \frac{1}{n^2}.
$$ 
A union bound now gives that whp the result holds for all $\log n \leq k \leq \frac{n}{3}$.

To prove (\ref{EAk1}), it suffices to show that $K_n$ contains at most $\frac{(n-1)! (\log n)^{k-1}}{(k-1)! 2^{k}}$ 2-factors with exactly $k$ cycles. We can count these as follows: Define an \textit{ordered 2-factor} to be a 2-factor together with an ordering of its cycles. We can count the number of ordered 2-factors by first choosing some $k$-tuple $(a_1, a_2, \ldots, a_k)$ and counting those ordered 2-factors whose cycles have lengths $a_1, a_2, \ldots, a_k$ in that order. This can be done by simply ordering $V(G)$ and placing vertices $1$ to $a_1$ in the first cycle, vertices $a_1 + 1$ to $a_1 + a_2$ in the second, etc. This procedure will count each ordered 2-factor of the appropriate type $(2a_1) (2a_2) \cdots (2a_k)$ times, and hence the number of these ordered 2-factors is $\frac{n!}{2^k a_1 a_2 \cdots a_k}$. Summing over all valid $k$-tuples, we have that the total number of ordered 2-factors with $k$ cycles is
$$\sum \frac{n!}{2^k} \prod_{i=1}^k \frac{1}{a_i} \leq \frac{n!}{2^k} \frac{k}{n} (\log n)^{k-1}$$
by Lemma ~\ref{fracsum}. But the number of ordered 2-factors (with $k$ cycles) is simply $k!$ times the total number of such 2-factors, and the result follows immediately. \endproof
We now need a lower bound on the total number of 2-factors. To do this we use the following well known result (see e.g.~the proof of Lemma 2 in~\cite{FK05}).
\begin{lemma}\label{perm}
Let $r$ be even, and let $H$ be an $r$-regular graph on $n$ vertices. Then $H$ contains at least $(\frac{r}{2n})^n n!$ $2$-factors.
\end{lemma}

\proof It is easy to see that the number of perfect matchings of a 
$d$-regular bipartite graph $B$  with vertex classes of size $n$ equals the permanent of the
incidence matrix of $B$. Egorychev~\cite{VDW1} and Falikman~\cite{VDW2} proved that the value of this permanent is at least
$\left(\frac{d}{n}\right)^n n!$ (this confirmed a conjecture of van der Waerden). So we take an orientation of $H$ in which every vertex has in- and out-degree $\frac{r}{2}$.\COMMENT{We can construct such an orientation by separating $H$ into 2-factors (with arbitrarily many cycles) and orienting each 2-factor individually.} Form a bipartite graph $B$ whose vertex classes $X$, $Y$ are each copies of $V(H)$, with an edge $xy$ for each $x \in X$, $y \in Y$ such that $\vec{xy}$ is an edge of the orientation of $H$. Now $B$ is $\frac{r}{2}$-regular and hence has at least $\left(\frac{r}{2n}\right)^n n!$ perfect matchings. But any perfect matching in $B$ yields a 2-factor in $H$, and distinct matchings yield distinct 2-factors. \endproof

\noindent \textit{Proof of Lemma ~\ref{2factor}.} It suffices to show that whp $A_{k_0+1} < (\frac{r}{2n})^n n!$, and hence by Lemma ~\ref{num2factors} it suffices that
$$\left(\frac{r}{2n}\right)^n n! \geq \frac{n! (\log n)^{2 k_0} p_0^n}{k_0! 2^{k_0}},$$
which holds as long as
$$\left(\frac{2np_0}{r}\right)^n \leq \frac{2^{k_0} k_0!}{(\log n)^{2 k_0}}.$$
Noting that $k_0! \geq (\frac{k_0}{e})^{k_0}$, it suffices that
\begin{align*}
n \log \frac{2np_0}{r} &\leq k_0(\log k_0 + \log 2 - 2 \log \log n -1) \\
&= \frac{\kappa n}{\log n}(\log n - 3 \log \log n + \log \kappa + \log 2 - 1).
\end{align*}
Since $\log \kappa + \log 2 - 1 > 0$, this follows immediately from
$$\log \frac{2np_0}{r} \leq \frac{\kappa}{2} \leq \kappa\left(1 - \frac{3 \log \log n}{\log n}\right).$$
\endproof
\begin{cor}\label{2fdecomp}
Let $m = \frac{1}{2}(1-\eta)np_0$. Then $G_1$ contains a collection of at least $m$ edge-disjoint $2$-factors, each with at most $k_0$ cycles.
\end{cor}
\proof By Corollary ~\ref{regular}, $G_1$ contains a regular subgraph $H$ of degree $r_1 = \left(1 - \frac{3\eta}{4}\right)np_0$. By Lemma ~\ref{2factor} (noting that $H$ is also a regular subgraph of $G_0$), we can remove 2-factors with at most $k_0$ cycles from $H$ one by one as long as the degree of the resulting graph remains above $2n p_0 e^{-\frac{\kappa}{2}}$. Recalling by (\ref{kappa}) that $e^{-\frac{\kappa}{2}} = \frac{\eta}{16}$, this gives us a collection of $\frac{1}{2}\left(r_1 - \frac{\eta n p_0}{8}\right) \ge \frac{1}{2}(1-\eta)np_0 = m$ 2-factors, each with at most $k_0$ cycles. \endproof


\section{Converting 2-factors into Hamilton cycles}
Applying Corollary ~\ref{2fdecomp} yields a collection $F_1, F_2, \ldots, F_m$ of edge-disjoint 2-factors in $G_1$, each with at most $k_0$ cycles, where 
\begin{equation}
m = \frac{1}{2}(1-\eta)np_0 \quad \textrm{ and } \quad k_0 = \frac{\kappa n}{\log n}.
\label{defnmk0}
\end{equation} 
Now we wish to convert these 2-factors into Hamilton cycles. Our proof develops ideas from Krivelevich and Sudakov \cite{mk-bs-2003}. Our strategy will be to show that for each $F_i$, we can connect the cycles of $F_i$ into a Hamilton cycle using edges of $G_0 \backslash (F_1 \cup F_2 \cup \ldots \cup F_m)$. We do this by incorporating the cycles of $F_i$ one by one into a long path, and then finally closing this path to a Hamilton cycle.

Let
$$E_0 = \frac{\log n}{\log(\frac{\eta w_0}{20})} \quad \textrm{ and } \quad E_1 = \frac{\log n}{\log(\frac{\eta^2 w_0}{10^5})}.$$

\begin{defin}\label{rotation}
Let $P = v_1 v_2 \ldots v_\ell$ be a path in $G_0$ with endpoints $v_1 = x, v_\ell = y$, and let $\Gamma$ be a spanning subgraph of $G_0$, whose edges are disjoint from those of $P$. Let $v_i$ be a vertex of $P$ such that $v_i y$ is an edge of $\Gamma$. A \textit{rotation} of $P$ about $y$ with \textit{pivot} $v_i$ is the operation of deleting the edge $v_i v_{i+1}$ from $P$ and adding the edge $v_i y$ to form a new path $v_1 v_2 \ldots v_i y v_{\ell-1} v_{\ell-2} \ldots v_{i+1}$ with endpoints $x$ and $y' = v_{i+1}$. Call the edge $v_i v_{i+1}$ the broken edge of the rotation. 
\end{defin}

Call a spanning subgraph $F$ of $G_0$ a \textit{broken $2$-factor} if $F$ consists of a collection of vertex-disjoint cycles together with a vertex-disjoint path, which we call the \textit{long path} of $F$. The key to our proof of Theorem ~\ref{AppHamDecomp} is the following lemma:

\begin{lemma}\label{rotext}
Let $P$ be a path in $G_0$. Let $\Gamma$ be a spanning subgraph of $G_0$ whose edges are disjoint from those of $P$, such that 
\begin{equation}
e(G_2 \backslash \Gamma) \leq \frac{\eta^6 n^2 p_0}{10^{17}} \label{remedges}
\end{equation}
and
\begin{equation}
\delta(\Gamma) \geq \frac{n \eta p_0}{4} \label{gammamindeg}.
\end{equation}
Then there exists a sequence of at most $2E_0 + 2E_1$ rotations which can be performed on $P$, using edges of $\Gamma$, to produce a new path $P'$, such that at least one of the following holds:
\begin{itemize}
\item[{\rm (i)}] $|P| \geq \frac{\eta n}{200}$ and the endpoints $x, y$ of $P'$ are joined by an edge of $\Gamma$.
\item[{\rm (ii)}] One of the endpoints $x, y$ of $P'$ is joined to a vertex outside $P'$ by an edge of~$\Gamma$.
\end{itemize}
\end{lemma}

Before we prove Lemma \ref{rotext}, we will first show how it is used to prove Theorem ~\ref{AppHamDecomp}. 
Our aim will be to convert each 2-factor $F_i$ into a Hamilton cycle $H_i$ in turn. 

For this, we using the following algorithm: Let $F^*$ be a broken 2-factor formed by removing an edge of $F_i$ arbitrarily. 
Let $j$ be the number of \textit{steps} performed so far during the conversion process (both on $F^*$ and on the $2$-factors which have already been converted into Hamilton cycles). Here a step is taken to mean a single application of Lemma ~\ref{rotext} to either obtain a broken $2$-factor with fewer cycles or to close a Hamilton path to a cycle. Let 
$$\Gamma_j = G_0 \backslash (H_1 \cup  \ldots \cup H_{i-1} \cup F^* \cup F_{i+1} \cup \ldots \cup F_m).$$ 
Then since $\delta(G_0) \geq \delta(G_1) \geq (1-\frac{\eta}{2})np_0$ by Lemma ~\ref{mindegree}, and 
$$\Delta(H_1 \cup  \ldots \cup H_{i-1} \cup F^* \cup F_{i+1} \cup \ldots \cup F_m) \leq 2m = (1-\eta)np_0,$$ 
we have that $\delta(\Gamma_j) \geq \frac{\eta n p_0}{4}$. Assume also that $e(G_2 \backslash \Gamma_j) \leq 4jE_1$, and that 
\begin{equation} \label{jbound}
j \leq k_0 m \stackrel{(\ref{defnmk0})}{\leq} \frac{\kappa n^2 p_0}{2 \log n}.
\end{equation} 
Then
\begin{equation} \label{G2gamma}
e(G_2 \backslash \Gamma_j) \stackrel{(\ref{jbound})}{\leq} \frac{2\kappa n^2 p_0}{\log(\frac{\eta^2 w_0}{10^5})}
\stackrel{(\ref{kappa})}{\leq} \frac{\eta^6 n^2 p_0}{10^{17}}.
\end{equation} 

Let $P^*$ be the long path of $F^*$. If $P^*$ is a Hamilton path, then Lemma ~\ref{rotext} applied with $\Gamma = \Gamma_j$ and $P = P^*$ shows that after at most $2E_0 + 2E_1$ rotations we can close $P^*$ to a Hamilton cycle $H_i$. We then move on to the next 2-factor $F_{i+1}$. If there are no 2-factors remaining (i.e., if $i=m$), then we have constructed the required set of $m$ edge-disjoint Hamilton cycles.

Otherwise by Lemma ~\ref{rotext}, after at most $2E_0 + 2E_1$ rotations we can either join an endpoint of $P^*$ to a vertex $x$ outside $P^*$, or we can close $P^*$ to form a cycle $C^*$. In the first case $x$ will be a vertex of some cycle $C_x$ of $F^*$, and we can delete one of the edges of $C_x$ incident to $x$ to form a new path $P^{**}$ which incorporates $C_x$. We then redefine $F^*$ to be the union of $P^{**}$ with the remaining cycles of $F_i$; this is a broken 2-factor with long path $P^{**}$, which has one cycle fewer than before. The algorithm then proceeds to the next step.

In the second case we have that $|C^*| = |P^*| \geq \frac{\eta n}{200}$. Now if $|C^*| \leq n - \frac{\eta n}{200}$, then by Lemma ~\ref{ABlinear} we have $e_{G_2}(V(C^*), [n] \backslash V(C^*)) \geq \frac{\eta^2 n^2 p_0}{50000}$. Since 
$$
e(G_2 \backslash \Gamma_j) 
\stackrel{(\ref{G2gamma})}{\leq} \frac{\eta^6 n^2 p_0}{10^{17}} 
< \frac{\eta^2 n^2 p_0}{50000},
$$ 
there must exist an edge in $\Gamma_j$ from some vertex $y$ of $C^*$ to a vertex outside $C^*$. On the other hand, if $|C^*| \geq n - \frac{\eta n}{200}$ then applying Lemma ~\ref{boundary} with $H = \Gamma_j$ and $A = [n] \backslash V(C^*)$ implies the same.\COMMENT{We have $\delta_A \geq \frac{\eta n p_0}{4} \geq 12 \log n$. So if $\frac{\log n}{ap_0} \geq \frac{7}{2}$, then $b=|B_\Gamma(A)| \geq a > 0$, and if $\frac{\log n}{ap_0} \leq \frac{7}{2}$, then $3a + b \geq \frac{\eta n}{28}$ and so $b \geq \frac{\eta n}{200} > 0$.} We then delete one of the edges of $C^*$ incident to $y$ and extend the resulting path as in the first case.

We run this algorithm until the last 2-factor $F_m$ has been converted into a Hamilton cycle. Now since each step either reduces the number of cycles in a broken $2$-factor or closes a Hamilton path to a Hamilton cycle, the algorithm will terminate after at most $k_0m$
steps. It remains to justify our assumption that $e(G_2 \backslash \Gamma_j) \leq 4j E_1$, for each $j$ (i.e., at each step). We can prove this by induction: $G_2 \subseteq \Gamma_0$, and since at most $2 E_0 + 2 E_1$ rotations are performed at each step, it follows that $e(\Gamma_j \backslash \Gamma_{j+1}) \leq 2 E_0 + 2 E_1 + 2\leq 4 E_1$.%
   \COMMENT{If we close the long path $P^*$ into a cycle $C^*$, we need one more edge for this and then one further edge to connect
$C^*$ to another cycle in the 2-factor}
So $e(G_2 \backslash \Gamma_{j+1}) \leq 4j E_1 + 4 E_1 = 4 (j+1)E_1$, as required. \endproof

It remains to prove Lemma \ref{rotext}. Our strategy will be as follows: We can assume that whenever we have an endpoint $x$ of a path $P'$ obtainable by fewer than $2E_0 + 2E_1$ rotations of $P$, then all of its neighbours lie on $P'$ (otherwise ~(ii) holds). So assuming this, we try to form some large sets $A$, $B$, such that for any $a \in A, b \in B$, we can obtain a path $P'$ with endpoints $a, b$. Then Lemma ~\ref{ABlinear} together with~(\ref{remedges}) will allow us to close $P'$ to a cycle.

We will (eventually) obtain the sets $A, B$ by dividing the path $P$ into two segments, and showing that we can perform a large number of rotations using only those pivots which lie all in one half or all in the other. This will ensure that the rotations involving the first endpoint do not interfere with those involving the second endpoint, and vice versa. In order to do this we need to show two things: Firstly, there exists a subset $C_1$ of the first segment of the path and a subset $C_2$ of the second segment of the path, such that for $i = 1, 2$, each vertex in $C_i$ has many neighbours which also lie in $C_i$. In fact since we are concerned with the successors or predecessors of the neighbours rather than the neighbours themselves, we will require the neighbours to lie in the interior (taken along $P$) of $C_i$. Secondly we will show that we can force the endpoints of the path to actually lie in these subsets.

We can accomplish the latter property by showing that the subsets are sufficiently large and by performing rotations until each endpoint lies in its corresponding subset. The obvious problem with this is that as we perform these rotations, $C_1$ and $C_2$ will cease to lie in their respective segments. So instead of defining $C_1, C_2$ immediately, we construct a subset $C$ of $V(P)$ with certain properties; then after rotating so that $a, b$ lie in $C$, we will define $C_1, C_2$ to be subsets of $C$, and the properties of $C$ will ensure that the vertices of each $C_i$ have many neighbours in $int(C_i)$. Here the \textit{interior} $int(C_i)$ of $C_i$ is the set of elements $x$ of $C_i$ such that both of the vertices adjacent to $x$ along $P$ also lie in $C_i$.

We start with the following lemma, where $k  = \log n$.

\begin{lemma}\label{subsetC}
Let $\varepsilon = \frac{\eta}{600}$, and $P \subseteq G_0$ be a path, $n' := |P| \geq \frac{\eta n}{200}$. Let $\Gamma$ be a spanning subgraph of $G_0$, edge-disjoint from $P$, which satisfies (\ref{remedges}). Let $W_1, W_2,\dots, W_k$ be a partition of $P$ into segments whose lengths are as equal as possible. Then there exists $S \subseteq [k]$ with $|S| \geq (1-\varepsilon)k$, and subsets $W'_i \subseteq W_i$ for each $i \in S$ with $|int(W'_i)| \geq (1-\varepsilon)\frac{n'}{k}$, such that for any $x \in W'_i$, and for at least $|S|-\varepsilon k$ of the sets $W'_j$, $|N_{\Gamma}(x) \cap int(W'_j)| \geq \frac{\eta p_0 n'}{20k}$.
\end{lemma}

\proof We start with $S = [k]$ and $W'_i = W_i$, and as long as there exists $i \in S$ and a vertex $x \in W'_i$, such that $|N_{\Gamma}(x) \cap int(W'_j)| \leq \frac{\eta p_0 n'}{20k}$ for at least $\varepsilon k$ values of $j \in S$, we remove $x$. (In this case, call $x$ \textit{weakly connected} to $W'_j$.) Further, if at any stage there exists $i \in S$ such that $|int(W'_i)| \leq (1-\varepsilon)\frac{n'}{k}$, then we remove $i$ from $S$. 

We claim that this process must terminate before $\frac{\varepsilon^2 n'}{4}$ vertices are removed. Indeed, suppose we have removed $\frac{\varepsilon^2 n'}{4}$ vertices and let $R$ be the set of removed vertices. Now $|R| = \frac{\varepsilon^2 n'}{4}$, and so $\sum_{i=1}^k |int(W'_i)| \geq (1-\frac{3\varepsilon^2}{4})n'$. Hence we have $|int(W'_i)| \geq (1-\varepsilon)\frac{n'}{k}$ for at least $1-\frac{3 \varepsilon k}{4}$ values of $i$\COMMENT{If there exist $k' > \frac{3 \varepsilon k}{4}$ values of $i$ such that $|int(W'_i)| < (1-\varepsilon)\frac{n'}{k}$, then $\sum_{i=1}^k |int(W'_i)| = \sum_{|int(W'_i)| < (1-\varepsilon)\frac{n'}{k}} |int(W'_i)| +  \sum_{|int(W'_i)| \geq (1-\varepsilon)\frac{n'}{k}} |int(W'_i)| < k'(1-\varepsilon)\frac{n'}{k} + (k-k')\frac{n'}{k} = n' - k' \varepsilon \frac{n'}{k} \leq (1-\frac{3\varepsilon^2}{4})n'$.}, i.e., at most $\frac{3 \varepsilon k}{4}$ indices have been removed from our original set $S$. So each $x \in R$ is still weakly connected to at least $\frac{\varepsilon k}{4}$ sets $W'_i$ with $i \in S$. For each $i \in S$, let $WC(i)$ be the set of vertices $x \in R$ which are weakly connected to $W'_i$.

Now consider the set $S_0 = \{i \in S \mid |WC(i)| \geq \frac{\varepsilon^3 n'}{32}\}$. Note that if $i \in S_0$, then
$$\frac{|int(W'_i)|}{n} \frac{|WC(i)|}{n} np_2 \geq \frac{(1-\varepsilon)\varepsilon^3 (n')^2 \eta w_0 \log n}{128 k n^2} \geq \frac{\eta^6 w_0}{10^{16}} \geq 700.$$
So Lemma ~\ref{ABlinear} implies that the number of edges of $G_2$ between $int(W'_i)$ and $WC(i)$ is at least $$\frac{13}{14} \frac{\eta p_0 |WC(i)|(1-\varepsilon)n'}{4k} \geq \frac{\eta p_0 |WC(i)|n'}{10k}.$$
But by the definition of $WC(i)$, $\Gamma$ contains at most $\frac{\eta p_0 |WC(i)|n'}{20k}$ edges between $int(W'_i)$ and $WC(i)$, and hence $G_2 \backslash \Gamma$ contains at least this many edges between $int(W'_i)$ and $WC(i)$.

Observe that $\sum_{i \in S} |WC(i)| \geq |R|\frac{\varepsilon k}{4} = \frac{\varepsilon^3 n'k}{16}$, since each $x \in R$ is weakly connected to $W'_i$ for at least $\frac{\varepsilon k}{4}$ values of $i \in S$. But since $\sum_{i \in S \backslash S_0} |WC(i)| \leq \frac{\varepsilon^3 n'k}{32}$, we have that $\sum_{i \in S_0} |WC(i)| \geq \frac{\varepsilon^3 n'k}{32}$. Hence $G_2 \backslash \Gamma$ contains at least 
$$
\frac{\eta n' p_0}{20 k} \sum_{i \in S_0} |WC(i)| 
\geq \frac{\eta \varepsilon^3 (n')^2 p_0}{640} 
\geq \frac{\eta^3 \varepsilon^3 n^2 p_0}{64 \cdot 4 \cdot 10^5} 
\geq \frac{\eta^6 n^2 p_0}{10^{16}}
$$ edges, which would contradict (\ref{remedges}). This proves the claim, and now we consider the sets $W'_i$ as they are at the point at which the process terminates. It is immediate that for each $i \in S$, $W'_i$ satisfies the requirements of the lemma. But since we have removed at most $\frac{3 \varepsilon k}{4}$ indices from our original set $S$, we also have $|S| \geq (1-\varepsilon)k$. \endproof

Let $C = \bigcup_{i \in S} W'_i$ and note that $|C| \geq (1 - \varepsilon)^2 n'$. We now need to show that the set of vertices which we can make into endpoints of $P$ with relatively few rotations is of size at least $2 \varepsilon n$. Doing this gives immediately that one of these endpoints must be an element of $C$.
\begin{lemma}\label{expansion}
Let $P$ be a path in $G_0$ with endpoints $a, b$, and $\Gamma$ be a spanning subgraph of $G_0$, edge-disjoint from $P$, which satisfies (\ref{gammamindeg}). Let $S_t$ be the set of vertices $x \in P\backslash \{b\}$ such that a path $P'$ with endpoints $x, b$ can be obtained from $P$ by at most $t$ rotations. Then $|S_{t+1}| \geq \frac{1}{2}|B_\Gamma(S_t)| - |S_t|$.
\end{lemma}
\proof For a vertex $x \in P$, let $x^-, x^+$ be the predecessor and successor of $x$ along $P$, respectively. Let $T = \{x \in B_\Gamma(S_t) \mid x^-, x^+ \notin S_t\}$. If $x \in T$, then since neither $x$ nor any of its neighbours on~$P$ are in $S_t$, the neighbours of $x$ are preserved by every sequence of at most $t$ rotations of $P$; i.e., $x^+$ and $x^-$ are adjacent to $x$ along any path obtained from $P$ by at most $t$ rotations. It follows that one of $x^-, x^+$ must be in $S_{t+1}$. Indeed, starting from $P$, we can perform $t$ rotations to obtain a path with endpoints $z, b$, such that $zx$ is an edge of $\Gamma$. Now by one further rotation with pivot $x$ and broken edge either $xx^+$ or $xx^-$, we obtain a path whose endpoints are either $x^+, b$ or $x^-, b$.

Now let $T^+ = \{x^+ \mid x \in T, x^+ \in S_{t+1}\}$ and $T^- = \{x^- \mid x \in T, x^- \in S_{t+1}\}$. It follows from the above that either $|T^+| \geq \frac{|T|}{2}$ or $|T^-| \geq \frac{|T|}{2}$, and both of these are subsets of $S_{t+1}$. Hence $|S_{t+1}| \geq \frac{|T|}{2} \geq \frac{1}{2}(|B_\Gamma(S_t)| - 2|S_t|)$. \endproof

\begin{cor}\label{initrots}
Either $|S_{E_0}| \geq \frac{\eta n}{200}$, or some element of $S_{E_0}$ has a neighbour in $\Gamma$ lying outside $P$ (or both).
\end{cor}

\proof It suffices to show that as long as $|S_t| \leq \frac{\eta n}{200}$, and assuming no element of $S_t$ has a neighbour outside $P$, we have that $|S_{t+1}| \geq \min \{\frac{\eta w_0}{20}|S_t|, \frac{\eta n}{200}\}$. We apply Lemma ~\ref{boundary}, setting
$H = \Gamma$ and $A = S_t$. Now in the notation of Lemma ~\ref{boundary}, $\delta_A \geq \delta(H[V(P)]) \geq \frac{\eta np_0}{4}$ by (\ref{gammamindeg}), and so we have that either (i) $|B_\Gamma(S_t)| \geq (\frac{\eta w_0}{8} - 3)|S_t|$, or (ii) $3|S_t| + |B_\Gamma(S_t)| \geq \frac{\eta n}{28}$. If (i) holds then 
$$|S_{t+1}| \geq \frac{1}{2}|B_\Gamma(S_t)| - |S_t| \geq (\frac{\eta w_0}{20} + 1)|S_t| - |S_t| = \frac{\eta w_0}{20}|S_t|.$$ On the other hand if (ii) holds then 
$$|S_{t+1}| \geq \frac{1}{2}|B_\Gamma(S_t)| - |S_t| \geq \frac{\eta n}{56} - \frac{5}{2}|S_t| \geq \frac{\eta n}{200}.$$ \endproof

Corollary~\ref{initrots} implies that if our path~$P$ in Lemma~\ref{rotext} satisfies $|P|<\frac{\eta n}{200}$, then
alternative~(ii) of Lemma~\ref{rotext} holds. So suppose that $|P|\ge \frac{\eta n}{200}$. Then we can apply Lemma~\ref{subsetC}
to obtain a set $C=\bigcup_{i \in S} W'_i$.
Now since $\frac{\eta n}{200} + |C| > n'=|P|$, we have that either alternative~(ii) of Lemma~\ref{rotext} holds, or we can obtain in at most $E_0$ rotations a path with endpoints $a', b$ such that $a' \in C$. Suppose we are in the latter case. Repeating the argument for $b$ gives us a path $P'''$ with endpoints $a', b' \in C$ which is obtained from $P$ by at most $2E_0$ rotations.

Call a segment $W_i$ of $P$ \textit{unbroken} if none of the rotations by which $P'''$ is obtained had their pivot in $W_i$. Note that each unbroken segment is still a segment of $P'''$ in the sense that the vertices are consecutive and their adjacencies along the path are preserved. Since we have arrived at the path $P'''$ by at most $2 E_0$ rotations, there are at least $k - 2 E_0$ unbroken segments $W_i$, and for at least $k - 2E_0 - \varepsilon k$ of these we have that $i \in S$. Noting that $E_0 \leq \frac{k}{10}$, we are still left with at least $\frac{3k}{5}$ unbroken segments $W_i$ for which $i \in S$. Let us relabel these segments $W_i$ according to their order along $P'''$, and take $C_1 = \bigcup_{i \leq \frac{3k}{10}} W'_i$ and $C_2 = \bigcup_{i > \frac{3k}{10}} W'_i$. Note that for any $x \in C$ (and in particular for $x \in C_1$ and for $a'$),
\begin{equation}
|N_{\Gamma}(x) \cap int(C_1)| \geq \frac{\eta p_0 n'}{20k}\left(\frac{3k}{10} - \varepsilon k\right) \geq \frac{\eta p_0 n'}{70} \geq \frac{\eta^2 n p_0}{14000}.
\label{nhoodint}
\end{equation}

Now let $x_0$ be a vertex separating $C_1, C_2$ along $P'''$, and let $x_0$ divide $P'''$ into paths $P_{a'}, P_{b'}$. Let $U_t$ be the set of endpoints of paths obtainable by at most $t$ rotations about $a'$ with pivots lying only in $int(C_1)$. So these rotations affect only $P_{a'}$, and $P_{b'}$ is left intact in each of the resulting paths. 

\begin{lemma}\label{endrots}
Suppose that $|U_t| \leq \frac{\eta^2 n}{10^6}$. Then 
$$|B_\Gamma(U_t) \cap int(C_1)| \geq \min \left\{\frac{\eta^2 n}{150000}, \frac{\eta^2 w_0}{40000}|U_t| \right\}.$$
\end{lemma}

\proof Let $u = |U_t|$ and $u' = |B_\Gamma(U_t) \cap int(C_1)|$. Consider the case $\frac{\log n}{up_0} \geq \frac{7}{2}$. Then similarly to the proof of Lemma ~\ref{boundary},
\begin{align*}
\frac{u \eta^2 n p_0}{14000} &\stackrel{(\ref{nhoodint})}{\leq} \sum_{x \in U_t} |N_\Gamma(x) \cap int(C_1)| \leq e_\Gamma(U_t, B_\Gamma(U_t) \cap int(C_1)) + 2 e_\Gamma(U_t) \\
&\leq 2(u + u')\log n + 4u \log n,
\end{align*}
whence 
$$u' \geq \frac{(\eta^2 n p_0 - 84000 \log n)u}{28000 \log n} 
= \frac{(\eta^2 w_0 - 84000)u}{28000} \geq \frac{\eta^2 w_0}{40000}|U_t|.$$
On the other hand, if $\frac{\log n}{up_0} \leq \frac{7}{2}$ then
\begin{align*}
\frac{u \eta^2 n p_0}{14000} &\leq \sum_{x \in U_t} |N_\Gamma(x) \cap int(C_1)| \leq e_\Gamma(U_t, B_\Gamma(U_t) \cap int(C_1)) + 2 e_\Gamma(U_t) \\
&\leq 7up_0(u+u') + 14u^2 p_0
\end{align*}
and so $3u + u' \geq \frac{\eta^2 n}{98000}$. Hence $u' \geq \frac{\eta^2 n}{150000}$. \endproof

\begin{cor}\label{endrots2}
$|U_{E_1}| \geq \frac{\eta^2 n}{10^6}.$
\end{cor}

\proof It suffices to prove that for each $t$ such that $|U_t| \leq \frac{\eta^2 n}{10^6}$, either $|U_{t+1}| \geq \frac{\eta^2 n}{10^6}$ or $|U_{t+1}| \geq \frac{\eta^2 w_0}{10^5}|U_t|$ (or both). Similarly to Lemma ~\ref{expansion}, we have that $|U_{t+1}| \geq \frac{1}{2}|B_\Gamma(U_t) \cap int(C_1)| - |U_t|$, and now Lemma ~\ref{endrots} immediately gives the result. \endproof

\noindent\textit{Proof of Lemma ~\ref{rotext}.}
Suppose first that $|P| \leq \frac{\eta n}{200}$. Then Corollary ~\ref{initrots} immediately implies that we can obtain a path, one of whose endpoints has a neighbour in $\Gamma$ lying outside $P$, in at most $E_0$ rotations. So we may assume that $|P| \geq \frac{\eta n}{200}$, and hence the conditions of Lemma ~\ref{subsetC} are satisfied. Now we proceed as above to obtain a path $P''' = P_{a'} \cup P_{b'}$, with endpoints $a', b'$ and with sets $C_1, C_2$ satisfying (\ref{nhoodint}), such that $a' \in C_1 \subseteq P_{a'}$ and $b' \in C_2 \subseteq P_{b'}$.

Let $U = U_{E_1}$. Now similarly, we can rotate about $b'$ using only pivots in $C_2$, to obtain another set $U'$ of endpoints in another $E_1$ rotations, such that $|U'| \geq \frac{\eta^2 n}{10^6}$. Now by Lemma ~\ref{ABlinear}, there are at least $\frac{13}{14} \frac{\eta^4 n^2 p_0}{10^{12}}$ edges of $G_2$ between $U$ and $U'$, and since by (\ref{remedges}) $G_2 \backslash \Gamma$ contains fewer edges than this, it follows that there exists an edge $xy$ of $\Gamma$ between $U$ and $U'$. Now by the definition of $U$ and $U'$, we can obtain a path $P''$ with endpoints $x, b'$ from $P'''$ by a sequence of at most $E_1$ rotations, none of which affect the second half $P_{b'}$ of the path $P'''$. From $P''$ we can obtain a path $P'$ with endpoints $x, y$ by at most $E_1$ rotations. \endproof

\section{Acknowledgements}
We thank the referees for their detailed comments.

\medskip

{\footnotesize \obeylines \parindent=0pt

Fiachra Knox, Daniela K\"{u}hn \& Deryk Osthus
School of Mathematics
University of Birmingham
Edgbaston
Birmingham
B15 2TT
UK

\begin{flushleft}
{\tt{E-mail addresses}:
\{knoxf,kuehn,osthus\}@maths.bham.ac.uk}
\end{flushleft} }

\end{document}